# Terminal Control Area Capacity Estimation Model Incorporating Structural Space

**Jeong Woo Park[1], Huiyang Kim[2,*]**

[1] Department of Advanced Air Transportation, Korea Aerospace University, Goyang, Republic of Korea

[2] School of Air Transport, Transportation and Logistics, Korea Aerospace University, Goyang, Republic of Korea

[*] Corresponding author: igreenbee@kau.ac.kr; 76-10, Hanggongdaehak-ro, Deogyang-gu, Goyang-si, Gyeonggi-do, Republic of Korea

**Declaration of Interest**

The authors declare that they have no known competing financial interests or personal relationships that could have appeared to influence the work reported in this paper.

**Acknowledgement**

This work was supported by the GRRC program of Gyeonggi province (GRRCAerospace2023-B03, Fusion Technology Research Center for Advanced Air Mobility).

# Terminal Control Area Capacity Estimation Model Incorporating Structural Space

**Abstract**

The continuous growth in global air traffic demand highlights the need to accurately estimate airspace capacity for efficiently using limited resources in air traffic management (ATM) systems. Although previous studies focused on either sector capacity based on air traffic controllers (ATCo) workload or runway throughput, studies on the unique structural and functional characteristics of terminal control area (TMA) remain lacking. In this study, capacity is defined as the maximum occupancy count. Further, a TMA capacity estimation model grounded in structural space conceptually defined as the space formed by instrument flight procedures and traffic characteristics is developed. Capacity is estimated from the temporal flight distance, which represents the physical length of arrival paths converted to flight time, and the average time separation at the runway threshold considering traffic proportions and aircraft mix. The proposed model is applied to the Jeju International Airport TMA (RWY 07/25) using one year of ADS-B trajectory data. The estimated capacities are 9.3 (RWY 07) and 6.9 (RWY 25) aircraft, and the differences are attributed to the temporal flight distance. Sensitivity analysis shows that capacity is shaped by aircraft speed and air traffic control (ATC) separations, which implies that operational measures such as speed restrictions or adjusted separations effectively enhance capacity even within physically constrained TMA. The model offers a practical, transparent, and quantitative framework for TMA capacity assessment and operational design.

## 1. Introduction

Global air traffic demand continues to rise steadily despite the temporary decline attributed to external shocks such as the global economic crisis and outbreaks of pandemics such as Covid-19, and this trend is expected to persist in the future (European Organization for the Safety of Air Navigation, 2022; Baneshi et al., 2024). The supply capability of airports and airspace should be expanded to address this increase in demand. However, an imbalance between demand and supply has emerged because of the limitations of supply expansion, which lead to social costs such as delays (Xu et al., 2021).

The International Civil Aviation Organization (ICAO) and many countries worldwide have adopted the global air navigation plan (GANP) to effectively respond to this growing demand for air traffic. To this end, they established the development direction for air traffic management (ATM) systems while presenting technical and operational solutions for effectively utilizing limited resources and expanding supply capability (International Civil Aviation Organization, 2005; International Civil Aviation Organization, 2008; International Civil Aviation Organization, 2016a). Given this background, the accurate estimation and evaluation of the accommodation capacity of currently operating ATM systems is essential for their improvement.

The capacity of airports and airspace are essential elements of the ATM system, and they have been extensively investigated by utilizing diverse approaches that reflect changes in the environment and technology. Although various models have been proposed (Harris, 1973; Schmidt, 1976; Odoni and Simpson, 1979; Janić and Tošić, 1982; International Civil Aviation Organization, 1984; Yang and Kim, 1994; Horonjeff et al., 2010; Juričić et al., 2011; Neufville and Odoni, 2013; Welch, 2015; Hanson et al., 2016; Zhang et al., 2016; Han et al., 2022; Chae et al., 2023; Federal Aviation Administration, 2025a), further discussion is required based on the extant studies on airspace capacity. The ICAO approaches airspace capacity from the perspective of airspace sectors, presenting an integrated methodology without distinguishing between en-route and terminal control area (TMA) (International Civil Aviation Organization, 1984). Although most prior studies focused on airspace sectors and proposed corresponding models (Schmidt, 1976; International Civil Aviation Organization, 1984; Welch, 2015; Hanson et al., 2016; Federal Aviation Administration, 2025a), to the best of the authors' knowledge, research addressing TMA capacity remains limited. Some studies focused on airport capacity, wherein they recognized TMA capacity as dependent on runway capacity while disregarding the unique structure and independent role of TMA (Horonjeff et al., 2010; Neufville and Odoni, 2013).

Previous studies estimated the capacity from the perspective of air traffic controllers (ATCo) workload in airspace sectors (Juričić et al., 2011; Zhang et al., 2016) or presented the number of aircraft handled per hour at specific points as the capacity (Janić and Tošić, 1982; Yang and Kim, 1994). However, only a limited number of approaches address a spatial perspective based on the structure of the airspace. Odoni and Simpson (1979), Janić and Tošić (1982), and Yang and Kim (1994) classified TMA as an independent airspace unit and focused on its functional role and structural attributes. Although these studies provided a theoretical foundation for developing TMA capacity research, their practical application remained limited.

TMA is an airspace established near major airports at points where the flight paths of departure and arrival aircraft intersect or merge, and its nature is distinct from that of the general en-route airspace (Visser, 1991). Functionally, TMA not only connects airports with en-routes but also regulates traffic volume and absorbs a portion of delays, managing the flow of air traffic (Zhou et al., 2016). The functions and roles of TMA are determined not only by the size and extent of the physical space but also by the structural characteristics of networked instrument flight procedures (International Civil Aviation Organization, 2020; Federal Aviation Administration, 2022a; Federal Aviation Administration, 2022b; Ministry of Land, Infrastructure and Transport Republic of Korea, 2024). Instrument flight procedures are designed to minimize the ATCo intervention of for aircraft operating within the TMA and effectively support the objectives of air traffic service (ATS) and strategies of the GANP (McElhatton ET AL., 1997). As indicated by ICAO, methods similar to those used for airspace sectors comprising en-routes can be applied to assess the TMA capacity (International Civil Aviation Organization, 1984). The physical availability of TMA and its unique functions and roles arising from the form and arrangement of instrument flight procedures must be considered, and therefore, it needs to be distinguished from en-route airspace and defined based on its spatial and structural characteristics.

This study proposes a capacity model that considers the characteristics of TMA and indicates that further research on TMA capacity and an approach from a spatial perspective is required. Unlike sector capacity models based on ATCo workload, the proposed model is defined based on structural space, such as the size of the airspace and instrument flight procedures. Further, it reflects traffic characteristics such as traffic proportions, aircraft type mix, and aircraft speed. This study presents a capacity model to comprehensively reflect traffic characteristics based on the space and structure of the TMA, which supplements the limitations of previous research, provides a foundation for evaluating ATM system efficiency, and offers the potential for practical application in future airspace design and operational strategy development.

The remainder of this paper is organized as follows: Section 2 reviews the previous studies on TMA capacity and outlines the key implications. Section 3 defines fundamental concepts and assumptions and develops the proposed capacity estimation model. Section 4 applies the model to the Jeju International Airport TMA and presents the results, including a sensitivity analysis of the major variables. Finally, Section 5 summarizes the main findings and concludes the study.

## 2. Literature Review

Capacity is defined as the maximum number of aircraft or throughput that can be accommodated by an ATM system of airports or airspace per unit time (International Civil Aviation Organization, 2016b; International Civil Aviation Organization, 2018). Capacity is determined by the combined effects of various factors, which include not only the structure of airport facilities and airspace and their operating conditions but also the complexity of air traffic flow, ATCo workload, and meteorological conditions. ATS providers determine and publish capacity within the range that does not exceed the amount of traffic volume that can be handled safely by the ATM system (International Civil Aviation Organization, 2018).

Based on this general definition, capacity is defined differently for airports and airspace depending on the perspective and approach of the researcher. At airports, capacity is expressed as the number of aircraft movements that the airport can handle per unit time (Harris, 1973; Horonjeff et al., 2010; Neufville and Odoni, 2013), whereas airspace capacity is defined either as the maximum number of aircraft that can enter the airspace per unit time (Schmidt, 1976; Odoni and Simpson, 1979; Janić and Tošić, 1982; Yang and Kim, 1994; Juričić et al., 2011; Welch, 2015; Zhang et al., 2016; International Civil Aviation Organization, 2018; Han et al., 2022; Chae et al., 2023) or as the maximum number of aircraft that can occupy the airspace during a specific time period (Hanson et al., 2016; International Civil Aviation Organization, 2018; Chae et al., 2023; Federal Aviation Administration, 2025a). This study describe the TMA capacity based on the concept of "maximum occupancy count" such that it reflects the physical and structural characteristics of the airspace and defines the capacity as the maximum number of aircraft that can occupy the structural space without ATCo intervention. Further, "structural space" is defined as a three-dimensional space structured around instrument flight procedures and air traffic characteristics. Although structures constituting TMA vary according to the spatial configuration of the space, they are important in that they directly and indirectly affect the traffic flow and capacity of the TMA. However, the structural space does not simply refer to the physical length of instrument flight procedures, and instead, it refers to the temporal conceptualization of length (hereinafter referred to as "temporal flight distance") where the physical length is converted into aircraft flight time. Further, the characteristics of air traffic also consider differences in the proportions of traffic volume through each entry point of TMA.

In the structural space, aircraft either traverse TMA via multiple entry points connected to the en-route and then approach the airport, or they depart the airport and enter the en-route (Visser, 1991; Vempati and Bonn, 2012; Zhou et al., 2016). Considering the traffic patterns of arrivals and departures, congestion and delay within TMA can be attributed to arrival aircraft (Janić and Tošić, 1982; Yang and Kim, 1994; Chae et al., 2023). Accordingly, this study incorporates only arrival flow into the model and excludes departures.

Research on airspace capacity, including that of TMA, has primarily focused on airspace sectors, recognizing ATCo workload as the main determinant (Schmidt, 1976; International Civil Aviation Organization, 1984; Juričić et al., 2011; Welch, 2015; Hanson et al., 2016; Zhang et al., 2016; Federal Aviation Administration, 2025a). In contrast, fewer studies treated TMA as an independent airspace and mathematically evaluated its capacity based on structural space (Janić and Tošić, 1982; Yang and Kim, 1994; Chae et al., 2023).

Janić and Tošić (1982) defined the flight routes of arrival aircraft based on standard instrument arrival routes (STARs) or actual trajectories and designated a point before reaching the runway threshold as the runway entry gate. They calculated the average time difference between successive aircraft at the runway entry gate. Their capacity model accounted for the geometric structure among flight routes, entry rates at different entry points, aircraft mix, and aircraft speeds along the different segments of the routes, in addition to other factors such as horizontal separation minima. The study applied the inverse of the calculated average time difference and presented the number of aircraft passing through the runway entry gate per hour as TMA capacity.

The model proposed by Janić and Tošić is significant in that it comprehensively considers both the structure of the airspace and traffic characteristics; however, its proposal to represent TMA capacity by the number of aircraft passing through a specific point can be interpreted as a throughput-based approach. Despite its academic contribution, the model has several limitations, including a constrained number of entry points and flight route segments, complexity related to calculating leading and trailing aircraft positions and inter-aircraft distances depending on route configurations, and the assumption of constant aircraft speeds along each segment, which makes it difficult to reflect the continuous changes in aircraft speed.

Building on the model of Janić and Tošić, Yang and Kim (1994) defined the TMA capacity as the number of arrival aircraft per hour at the runway threshold. Their study is similar to that of Janić and Tošić in that it is grounded in the concept of throughput and incorporates factors such as the structure of flight routes and aircraft mix. However, it uses a polar coordinate system to vectorize aircraft movements for simplifying mathematical complexity associated with the geometric structure of routes. Although the model simplifies the model of Janić and Tošić, it assumes that leading and trailing aircraft follow the same route and cannot realistically reflect variations in the speeds of arrival aircraft.

Chae et al. (2023) defined TMA capacity as the throughput and instantaneous capacity of arrival aircraft based on flight routes comprising STARs and instrument approach procedures (IAPs). In this context, throughput refers to the maximum number of arrival aircraft per hour at the runway threshold, while instantaneous capacity represents the maximum number of arrival aircraft occupying TMA at a given time. Chae et al. incorporated flight routes composed of STARs and IAPs, traffic characteristics, and other relevant factors to derive the average separation between aircraft and determine the throughput and instantaneous capacity on that basis. The average separation was expressed in units of time, with the assumptions that leading and trailing aircraft at the runway threshold satisfy the prescribed horizontal separation minima and that aircraft decelerate at a constant rate along routes, except within segments where speed restrictions are specified.

Chae et al.'s model not only considers the structure of TMA and continuous variations in the speeds of arrival aircraft but also defines TMA capacity in terms of both hourly throughput and instantaneous capacity from a spatial perspective. However, this model has several limitations in that it does not clearly define capacity as a single concept, fails to consider the effects of aircraft type by assuming identical average speeds for each flight route, and restricts some combinations of leading and trailing aircraft along the routes.

Juričić et al. (2011) defined TMA capacity using a simulation to estimate the traffic volume at which the ATCo workload reaches a certain level. Zhang et al. (2016) applied a mathematical model based on an ATCo workload to estimate the TMA capacity, defining the capacity as the number of aircraft corresponding to the maximum workload. Han et al. (2022) developed a regression model for predicting TMA capacity using machine learning algorithms (LightGBM and NGBoost) and presented the demand for arrivals and departures at the convergence point of the airport arrival rate and airport departure rate as the capacity.

Based on a review of previous studies, this study developed a TMA capacity model that considers both arrival paths from TMA entry points to runway threshold and arrival time difference between pairs of leading and trailing aircraft. This study is unlike the previous work because it focuses on developing a model that reflects the structural characteristics of the airspace and realistic traffic conditions while considering the following:

1) The TMA capacity is defined not in terms of throughput per unit time but as the maximum number of arrival aircraft that can occupy the structural space without considering ATCo workload.

2) Aircraft speeds are applied by type, with speeds set as a piecewise constant-rate deceleration from the entry point to the runway threshold, incorporating actual variations in aircraft speed during the arrival phase into the model.

3) All possible combinations of leading and trailing aircraft across different flight paths are considered to determine the effect of the geometric structure of the arrival flight path on capacity. In this process, longitudinal separation is applied to simplify the complexity of distance calculations based on aircraft positions (Brooker, 1983; Szurgyi et al., 2008).

Such a differentiated approach aligns with the study, which proposes capacity from the perspectives of structure and space and presents a mathematical model that enables its quantitative estimation, contributing to the efficient utilization of the ATM system and seamless flow of air traffic.

## 3. TMA Capacity Estimation Model

In this study, the framework and key factors required for capacity calculation were designed and specified in accordance with the defined concept of capacity to develop a TMA capacity model based on structural space. However, the ATCo workload was not considered to focus on the relationship between the TMA structure and capacity.

### 3.1 Basic Concept of Model Design

In accordance with the definitions of structural space and TMA capacity established in this study, the capacity model was designed based on the following fundamental concepts.

The arrival aircraft entering the TMA from multiple entry points follow their flight paths in compliance with the prescribed separation minima and terminate at the runway threshold (International Civil Aviation Organization, 2016b; Federal Aviation Administration, 2025b). Inter-aircraft separation can be expressed in either distance or time (Ministry of Land, Infrastructure and Transport Republic of Korea, 2022), and in this study, a time-based separation is applied. This choice is based on the consideration that aircraft speeds vary with position along the flight path and continuously decelerate as they approach the runway, making it difficult to apply a uniform distance-based separation.

However, the time separation between aircraft along the arrival flight path must meet or exceed the separation at the runway threshold to satisfy the condition that capacity estimation assumes no ATCo intervention such as radar vectoring. Therefore, this study applies the average time separation between aircraft derived at the arrival runway threshold as the separation between the aircraft along the arrival flight path (Fig. 1).

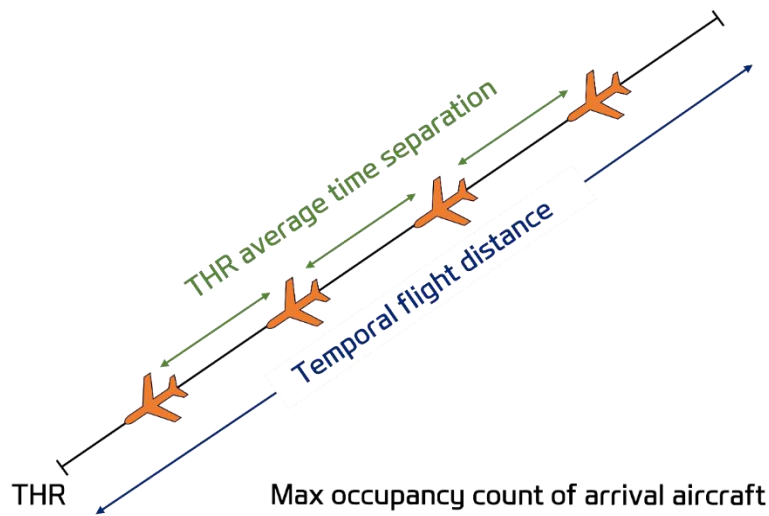

Fig. 1 Concept of terminal control area (TMA) capacity

The basic concept for the model design can be described as follows (Fig. 2):

(i) Identify key variables for capacity estimation and extract both traffic proportions at each entry point and speed data extracted by aircraft type.

(ii) Based on extracted data, convert the physical distance of the arrival flight path measured from each entry point to the runway threshold into the temporal flight distance.

(iii) Derive average time separation at the runway threshold by accounting for the combinations of leading and trailing aircraft using previously extracted traffic proportions and speed data.

(iv) Under the assumption that all aircraft are required to satisfy the average separation obtained at the runway threshold, determine the TMA capacity as the maximum number of arrival aircraft that can occupy the arrival flight path.

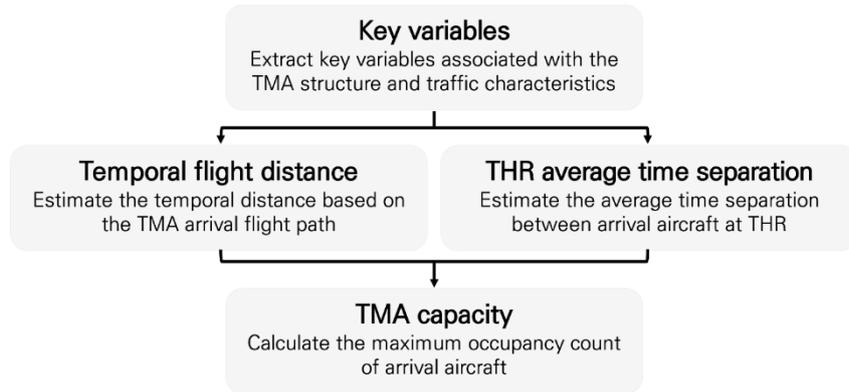

Fig. 2 TMA capacity estimation flowchart

3.2 Model Design Assumptions

This study makes the following assumptions for model design based on the basic concept:

(i) The arrival flight path is determined by ATS routes of the instrument flight procedures that connect each TMA entry point to the runway threshold.

(ii) Arrival aircraft continuously enter the TMA. From the entry point to the runway threshold, they follow only the determined ATS routes. They do not use any other routes and are not radar-vectored by the ATCo.

(iii) The aircraft speed profile is divided into two segments. The first segment extends from the entry point to the first merging point within the IAP where all paths converge (IAP merging point), and the second segment extends from that point to the runway threshold. It is assumed that aircraft decelerate at a constant rate defined for each segment.

(iv) Although at least one of lateral, longitudinal, or vertical separation may apply between aircraft (Brooker, 1983; Szurgyi et al., 2008), only longitudinal separation is considered in this model. For aircraft in close lateral proximity, the separation is assumed to be satisfied by vertical separation.

3.3 Model Design and Development

The number of arrival aircraft occupying TMA can be expressed as the count of aircraft that traverse the flight path from the entry point to the runway threshold while satisfying the prescribed inter-aircraft separation. According to the model design and its assumptions, the maximum occupancy count is determined by dividing the temporal flight distance by the average time separation at the runway threshold. TMA capacity ($\lambda_{\text{rwy}}$) expressed as the maximum occupancy count is defined as

$$\lambda_{\text{rwy}} = \frac{D_{\text{temp}}}{\bar{T}_{\text{thr}}}, \tag{1}$$

where $D_{\text{temp}}$ and $\bar{T}_{\text{thr}}$ represent the temporal flight distance and average time separation at the runway threshold, respectively.

The temporal flight distance is derived from the physical distance of flight paths comprising the TMA structure while accounting for both flight time and traffic proportion. Even when flight paths have the same physical distance, differences in the flight time and traffic proportion of those paths cause variations in the temporal flight distance. The arrival flight path within the TMA includes multiple paths, and the temporal flight distance for the overall arrival flight path is expressed as the distance obtained by summing each path. In this context, the temporal flight distance obtained by summing the individual arrival flight path can be represented as a single extended flight path. The temporal flight distance is given by

$$D_{\text{temp}} = \sum_{r=1}^{R}(\rho_r \cdot \bar{t}_{\text{tot},1}) + \sum_{r=2}^{R}\sum_{i=r}^{R}\rho_i \cdot (\bar{t}_{\text{tot},r} - \bar{t}_{\text{tot},r-1}), \tag{2}$$

where $R$, $r$, $\rho_r$, and $\bar{t}_{\text{tot},r}$ represent the number of arrival flight paths connected to the runway threshold, index of the arrival flight path when the paths are sorted in an ascending order by average flight time, traffic proportion of arrival flight path $r$, and average flight time of arrival flight path $r$, respectively.

The flight time for each arrival flight path varies based on the physical distance of the path and speed for each aircraft type. The average flight time of an arrival flight path ($\bar{t}_{\text{tot}}$) is expressed by accounting for the aircraft mix proportions and corresponding flight times for that path as

$$\bar{t}_{\text{tot}} = \sum_{c}(p_c \cdot t_{\text{tot},c}), \tag{3}$$

where $c$, $p_c$, and $t_{\text{tot},c}$ represent the aircraft type using the arrival flight path, proportion of aircraft type $c$ on the arrival flight path, and flight time of aircraft type $c$ along the arrival flight path, respectively.

This study divides each arrival flight path into two segments according to the aircraft speed profile.

The flight time of an arrival flight path ($t_{\text{tot}}$) is represented as the sum of the segmental flight times calculated by considering the length and aircraft speed of each segment as

$$t_{\text{tot}} = \frac{d_{\text{E-MP}_{\text{iap}}}}{\left(\frac{v_{\text{E}}+v_{\text{MP}_{\text{iap}}}}{2}\right)} + \frac{d_{\text{MP}_{\text{iap}}\text{-THR}}}{\left(\frac{v_{\text{MP}_{\text{iap}}}+v_{\text{THR}}}{2}\right)}, \quad (4)$$

where $d_{\text{E-MP}_{\text{iap}}}$, $d_{\text{MP}_{\text{iap}}\text{-THR}}$, $v_{\text{E}}$, $v_{\text{MP}_{\text{iap}}}$, and $v_{\text{THR}}$ represent the length from the entry point to $\text{MP}_{\text{iap}}$, length from $\text{MP}_{\text{iap}}$ to the runway threshold, aircraft speed at the entry point, aircraft speed at $\text{MP}_{\text{iap}}$, and aircraft speed at the runway threshold, respectively.

The average time separation at the runway threshold ($\bar{T}_{\text{thr}}$) is determined by the difference between the threshold-crossing times of the leading and trailing aircraft along the arrival flight path. This difference depends on the selected combination of the arrival flight path and aircraft type for the aircraft pair. $\lambda_{\text{rwy}}$ is inversely proportional to $\bar{T}_{\text{thr}}$, and therefore, minimizing $\bar{T}_{\text{thr}}$ is necessary for determining the maximum occupancy count. Accordingly, $\bar{T}_{\text{thr}}$ is defined by incorporating both the proportion associated with the combination of the leading and trailing aircraft and corresponding runway threshold-crossing time difference as

$$\bar{T}_{\text{thr}} = \sum_{k,l}\sum_{i,j}(\rho_k \cdot \rho_l \cdot p_{i/k} \cdot p_{j/l}) \cdot \Delta T_{i/k,j/l}, \quad (5)$$

where $k$, $l$, $i$, $j$, $\rho_k$, $\rho_l$, $p_{i/k}$, $p_{j/l}$, and $\Delta T_{i/k,j/l}$ represent the arrival flight path used by the leading aircraft, arrival flight path used by the trailing aircraft, aircraft type of the leading aircraft, aircraft type of the trailing aircraft, traffic proportion of the arrival flight path $k$, traffic proportion of the arrival flight path $l$, proportion of aircraft type $i$ on the arrival flight path $k$, proportion of aircraft type $j$ on the arrival flight path $l$, and runway threshold-crossing time difference between the leading and trailing aircraft, respectively.

$\bar{T}_{\text{thr}}$ is proportional to the runway threshold-crossing time difference between the leading and trailing aircraft ($\Delta T_{i/k,j/l}$). Therefore, minimizing $\bar{T}_{\text{thr}}$ involves minimizing $\Delta T_{i/k,j/l}$. The leading–trailing relationship for an aircraft pair flying on $k$ and $l$ is determined at the point where their common path begins ($\text{MP}_{kl}$) because this model considers only longitudinal separation between aircraft. Before entering the common path, the trailing aircraft does not follow the leading aircraft, and therefore, the satisfaction of the longitudinal separation is not considered. Longitudinal separation is considered only within the common path. When $k$ and $l$ correspond to the same arrival flight path, $\mathrm{MP}_{kl}$ is defined as the entry point of that path. If $\mathrm{MP}_{\mathrm{iap}}$ is the starting point of the common path of $k$ and $l$, $\mathrm{MP}_{kl}$ is set equal to $\mathrm{MP}_{\mathrm{iap}}$.

The common path of arrival flight paths is divided with respect to $\mathrm{MP}_{\mathrm{iap}}$ into two common subpaths (Common subpaths 1 and 2), which extend from $\mathrm{MP}_{kl}$ to $\mathrm{MP}_{\mathrm{iap}}$ and from $\mathrm{MP}_{\mathrm{iap}}$ to the runway threshold, respectively. However, if $\mathrm{MP}_{kl} = \mathrm{MP}_{\mathrm{iap}}$, common subpath 1 does not exist, and the common path includes only common subpath 2. Along the common path, the trailing aircraft follows the leading aircraft; when the leading aircraft enters $\mathrm{MP}_{kl}$, the trailing aircraft is located upstream of that point. At that moment, the along-path distance from the position of the trailing aircraft to $\mathrm{MP}_{kl}$ is defined as the initial spacing distance, and the corresponding flight time is defined as the initial spacing time (Fig. 3).

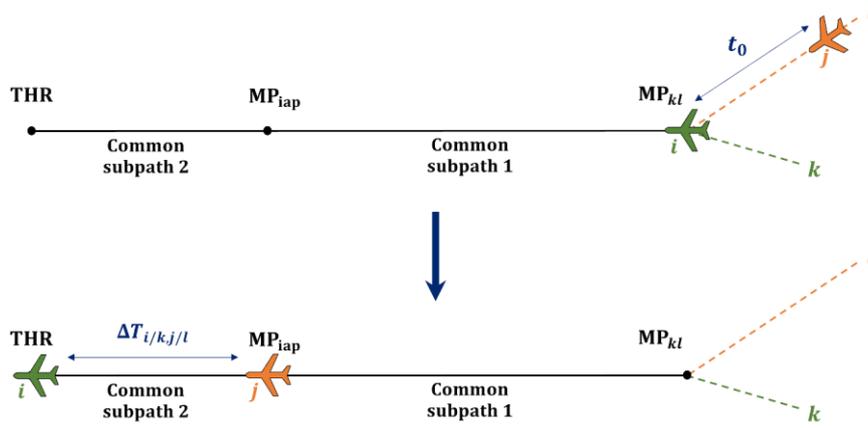

Fig. 3 Conceptual overview of $t_0$ and $\Delta T_{i/k,j/l}$

Accordingly, $\Delta T_{i/k,j/l}$ for a given leading–trailing aircraft combination is obtained by considering the flight time along the common path and initial spacing time as

$$\Delta T_{i/k,j/l} = ({}_{j/l}t_{\mathrm{com1}} + {}_{j/l}t_{\mathrm{com2}} + t_0) - ({}_{i/k}t_{\mathrm{com1}} + {}_{i/k}t_{\mathrm{com2}}), \tag{6}$$

where ${}_{i/k}t_{\mathrm{com1}}$, ${}_{i/k}t_{\mathrm{com2}}$, ${}_{j/l}t_{\mathrm{com1}}$, ${}_{j/l}t_{\mathrm{com2}}$, and $t_0$ represent the flight time of the leading aircraft along common subpath 1, flight time of the leading aircraft along common subpath 2, flight time of the trailing aircraft along common subpath 1, flight time of the trailing aircraft along common subpath 2, and initial spacing time between the leading and trailing aircraft with reference to $\text{MP}_{kl}$, respectively.

The flight times of the leading ($_{i/k}t_{\text{com1}}$, $_{i/k}t_{\text{com2}}$) and trailing aircraft ($_{j/l}t_{\text{com1}}$, $_{j/l}t_{\text{com2}}$) along the respective common subpaths are derived from the lengths of each subpath and the speeds of the leading and trailing aircraft. They are respectively expressed as

$$_{i/k}t_{\text{com1}} = \frac{_{kl}d_{\text{com1}}}{\left(\frac{_{i/k}v_{\text{MP}_{kl}} + _{i/k}v_{\text{MP}_{\text{iap}}}}{2}\right)}, \tag{7}$$

$$_{j/l}t_{\text{com1}} = \frac{_{kl}d_{\text{com1}}}{\left(\frac{_{j/l}v_{\text{MP}_{kl}} + _{j/l}v_{\text{MP}_{\text{iap}}}}{2}\right)}, \tag{8}$$

$$_{i/k}t_{\text{com2}} = \frac{_{kl}d_{\text{com2}}}{\left(\frac{_{i/k}v_{\text{MP}_{\text{iap}}} + _{i/k}v_{\text{THR}}}{2}\right)}, \tag{9}$$

$$_{j/l}t_{\text{com2}} = \frac{_{kl}d_{\text{com2}}}{\left(\frac{_{j/l}v_{\text{MP}_{\text{iap}}} + _{j/l}v_{\text{THR}}}{2}\right)}, \tag{10}$$

where $_{kl}d_{\text{com1}}$, $_{kl}d_{\text{com2}}$, $_{i/k}v_{\text{MP}_{kl}}$, $_{j/l}v_{\text{MP}_{kl}}$, $_{i/k}v_{\text{MP}_{\text{iap}}}$, $_{j/l}v_{\text{MP}_{\text{iap}}}$, $_{i/k}v_{\text{THR}}$, and $_{j/l}v_{\text{THR}}$ represent the length of common subpath 1, length of common subpath 2, speed of the leading aircraft at $\text{MP}_{kl}$, speed of the trailing aircraft at $\text{MP}_{kl}$, speed of the leading aircraft at $\text{MP}_{\text{iap}}$, speed of the trailing aircraft at $\text{MP}_{\text{iap}}$, speed of the leading aircraft at the runway threshold, and speed of the trailing aircraft at the runway threshold, respectively.

Aircraft flying along the arrival flight path decelerate uniformly according to their speed profiles. Thus, the speeds of the leading and trailing aircraft at $\text{MP}_{kl}$ ($_{i/k}v_{\text{MP}_{kl}}$, $_{j/l}v_{\text{MP}_{kl}}$) are determined based on the length along the arrival flight path from the entry point to $\text{MP}_{kl}$ by incorporating the entry-point speed of the aircraft into the path and its acceleration along the path.

$$_{i/k}v_{\text{MP}_{kl}} = \sqrt{_{i/k}v_{\text{E}}^2 - 2|_{i/k}a_{\text{E-MP}_{\text{iap}}}|\,_k d_{\text{E-MP}_{kl}}}, \tag{11}$$

$$_{j/l}v_{\text{MP}_{kl}} = \sqrt{_{j/l}v_{\text{E}}^2 - 2|_{j/l}a_{\text{E-MP}_{\text{iap}}}|\,_l d_{\text{E-MP}_{kl}}}, \tag{12}$$

where $_{i/k}v_{\text{E}}$, $_{j/l}v_{\text{E}}$, $_{i/k}a_{\text{E-MP}_{\text{iap}}}$, $_{j/l}a_{\text{E-MP}_{\text{iap}}}$, $_k d_{\text{E-MP}_{kl}}$, and $_l d_{\text{E-MP}_{kl}}$ represent the speed of the leading aircraft at the entry point, speed of the trailing aircraft at the entry point, acceleration of the leading aircraft over the segment of $k$ from the entry point to $\text{MP}_{\text{iap}}$, acceleration of the trailing aircraft over the segment of $l$ from the entry point to $\text{MP}_{\text{iap}}$, length of $k$ from the entry point to $\text{MP}_{kl}$, and length of $l$ from the entry point to $\text{MP}_{kl}$, respectively.

The accelerations of an aircraft for each segment of the arrival flight path divided according to the speed profile of the aircraft ($a_{\text{E-MP}_{\text{iap}}}$, $a_{\text{MP}_{\text{iap}}\text{-THR}}$) are calculated using segment lengths and aircraft speeds along the path, which are respectively expressed as

$$a_{\text{E-MP}_{\text{iap}}} = \frac{v_{\text{MP}_{\text{iap}}}^2 - v_{\text{E}}^2}{2 d_{\text{E-MP}_{\text{iap}}}}. \tag{13}$$

$$a_{\text{MP}_{\text{iap}}\text{-THR}} = \frac{v_{\text{THR}}^2 - v_{\text{MP}_{\text{iap}}}^2}{2 d_{\text{MP}_{\text{iap}}\text{-THR}}}. \tag{14}$$

$\Delta T_{i/k,j/l}$ is positively and linearly related to the initial spacing time ($t_0$), and therefore, minimizing $\Delta T_{i/k,j/l}$ requires that $t_0$ be minimized. The minimum $t_0$ is determined such that, while the leading aircraft traverses the common path of $k$ and $l$, the leading and trailing aircraft satisfy the air traffic control (ATC) separations applied in the TMA and at the runway threshold. The minimum $t_0$ is obtained as the smallest value that enables the aircraft pair to satisfy prescribed separations throughout the time interval $t \in [0, {}_{i/k}t_{\text{com1}} + {}_{i/k}t_{\text{com2}}]$ by setting the passage time of the leading aircraft at $\text{MP}_{kl}$ to $t = 0$ and its passage time at the runway threshold to $t = {}_{i/k}t_{\text{com1}} + {}_{i/k}t_{\text{com2}}$. This condition is formulated as the following optimization problem.

$$\min t_0$$

Subject to $\qquad\qquad\qquad\qquad\qquad\qquad\qquad\qquad\qquad\qquad\qquad\qquad\qquad\qquad\qquad (15)$

$$D_{i/k,j/l,n}(t) \geq S, \quad \forall n \in [1, \cdots, N], t \in [t_{n,\text{start}}, t_{n,\text{end}}]$$

$$D_{i/k,j/l,N}(t_{N,\text{end}}) \geq S_{\text{thr}}$$

where $S$, $S_{\text{thr}}$, $N$, $n$, $t_{n,\text{start}}$, $t_{n,\text{end}}$, $t_{N,\text{end}}$, and $D_{i/k,j/l,n}(t)$ represent the ATC separation applied in the TMA, ATC separation applied at the runway threshold, total number of subintervals into which the time interval $[0, {}_{i/k}t_{\text{com1}} + {}_{i/k}t_{\text{com2}}]$ is divided, index of the subinterval, start time of the $n$-th subinterval, end time of the $n$-th subinterval, end time of the last subinterval, and distance between the leading and trailing aircraft at time $t$ within the $n$-th subinterval, respectively.

Further, $S$ is determined according to the ATC separation applied in the target TMA, such as wake turbulence separation minima by aircraft type and radar horizontal separation minima (5 NM or 3 NM) (International Civil Aviation Organization, 2016b; Ministry of Land, Infrastructure and Transport Republic of Korea, 2022; Federal Aviation Administration, 2025b). $S_{\text{thr}}$ is determined by considering both $S$ and the characteristics of the runway. For example, in the case of a landing-only runway, $S_{\text{thr}}$ is generally set equal to $S$, whereas a larger value may be applied for a mixed-use runway handling both arrivals and departures.

Along the arrival flight path, the acceleration of the aircraft is piecewise, with a change at $\text{MP}_{\text{iap}}$. To reflect this characteristic, subintervals $[t_{n,\text{start}}, t_{n,\text{end}}]$ are defined by partitioning the time interval $[0, {}_{i/k}t_{\text{com1}} + {}_{i/k}t_{\text{com2}}]$, which is measured as the flight time of the leading aircraft along the common path with respect to the passage times of the leading and trailing aircraft at $\text{MP}_{\text{iap}}$. Based on $t_0$, $[0, {}_{i/k}t_{\text{com1}} + {}_{i/k}t_{\text{com2}}]$ is partitioned into one to three subintervals ($N = 1 - 3$), and they are expressed as

$$t \in \begin{cases} [0, {}_{i/k}t_{\text{com1}}] \cup [{}_{i/k}t_{\text{com1}}, {}_{i/k}t_{\text{com1}} + t_0] \cup [{}_{j/l}t_{\text{com1}} + t_0, {}_{i/k}t_{\text{com1}} + {}_{i/k}t_{\text{com2}}], & \text{if } {}_{i/k}t_{\text{com1}} + {}_{i/k}t_{\text{com2}} > {}_{j/l}t_{\text{com1}} + t_0 \\ [0, {}_{i/k}t_{\text{com1}}] \cup [{}_{i/k}t_{\text{com1}}, {}_{i/k}t_{\text{com1}} + {}_{i/k}t_{\text{com2}}], & \text{if } {}_{i/k}t_{\text{com1}} + {}_{i/k}t_{\text{com2}} \leq {}_{j/l}t_{\text{com1}} + t_0 \end{cases}, \quad (16)$$

where $t_0$ must satisfy the conditions

$$t_0 > 0,$$

$$t_0 > {}_{i/k}t_{\text{com1}} - {}_{j/l}t_{\text{com1}},$$

$$t_0 > ({}_{i/k}t_{\text{com1}} + {}_{i/k}t_{\text{com2}}) - ({}_{j/l}t_{\text{com1}} + {}_{j/l}t_{\text{com2}}),$$

$D_{i/k,j/l,n}(t)$ is the function of time $t$ that indicates the distance between the leading and trailing aircraft along the common path. This distance is formulated in terms of the initial distance for $[t_{n,\text{start}}, t_{n,\text{end}}]$, as well as the speeds and accelerations of the leading and trailing aircraft over that subinterval. It is defined as a quadratic function of $t$ and is expressed as

$$D_{i/k,j/l,n}(t) = s_n + ({}_{i/k}v_{n,\text{start}} - {}_{j/l}v_{n,\text{start}})(t - t_{n,\text{start}}) - \frac{1}{2}(|{}_{i/k}a_n| - |{}_{j/l}a_n|)(t - t_{n,\text{start}})^2, \quad (17)$$

where $s_n$, $_{i/k}a_n$, $_{j/l}a_n$, $_{i/k}v_{n,\text{start}}$, and $_{j/l}v_{n,\text{start}}$ represent the initial distance for the $n$-th subinterval between the leading and trailing aircraft, acceleration of the leading aircraft for the $n$-th subinterval, acceleration of the trailing aircraft for the $n$-th subinterval, speed of the leading aircraft at $t_{n,\text{start}}$, and speed of the trailing aircraft at $t_{n,\text{start}}$, respectively.

The initial distance for $[t_{n,\text{start}}, t_{n,\text{end}}]$ $s_n$ represents the distance between the leading and trailing aircraft along the common path at $t = t_{n,\text{start}}$. In the first subinterval ($n = 1$), the initial distance $s_1$ represents the initial spacing distance at $\text{MP}_{kl}$, which is determined by $t_0$. For subsequent subintervals ($n > 1$), the initial distance $s_n$ is derived through the function of the previous subinterval, $D_{i/k,j/l,n-1}(t)$, and the expression for $s_n$ is formulated as

$$s_n = \begin{cases} \left(_{j/l}v_{\text{MP}_{kl}} + \frac{1}{2}|_{j/l}a_1|t_0\right) \cdot t_0, & \text{if } n = 1 \\ D_{i/k,j/l,n-1}(t_{n,\text{start}}), & \text{if } n > 1 \end{cases}, \quad (18)$$

For $[t_{n,\text{start}}, t_{n,\text{end}}]$, the leading and trailing aircraft travel at distinct constant accelerations ($_{i/k}a_n$, $_{j/l}a_n$). $_{i/k}a_n$ and $_{j/l}a_n$ are determined as either $a_{\text{E-MP}_{\text{iap}}}$ or $a_{\text{MP}_{\text{iap}}\text{-THR}}$ according to the positions of the leading and trailing aircraft along the arrival flight path, as specified in Eq. (13) or (14). When $t_0$ causes the trailing aircraft to be positioned before the entry point, the acceleration of the trailing aircraft from that position to the entry point is set equal to $a_{\text{E-MP}_{\text{iap}}}$. The speeds of the leading and trailing aircraft at $t = t_{n,\text{start}}$ ($_{i/k}v_{n,\text{start}}$, $_{j/l}v_{n,\text{start}}$) are obtained from $t_0$ and the accelerations and lengths of the subinterval, given as

$$_{i/k}v_{n,\text{start}} = \begin{cases} _{i/k}v_{\text{MP}_{kl}} & , \text{if } n = 1 \\ _{i/k}v_{n-1,\text{start}} - |_{i/k}a_{n-1}| \cdot (t_{n-1,\text{end}} - t_{n-1,\text{start}}), & \text{if } n > 1 \end{cases}, \quad (19)$$

$$_{j/l}v_{n,\text{start}} = \begin{cases} _{j/l}v_{\text{MP}_{kl}} + |_{j/l}a_1|t_0 & , \text{if } n = 1 \\ _{j/l}v_{n-1,\text{start}} - |_{j/l}a_{n-1}| \cdot (t_{n-1,\text{end}} - t_{n-1,\text{start}}), & \text{if } n > 1 \end{cases}, \quad (20)$$

Fundamentally, the optimization problem for obtaining the minimum $t_0$ is to find the smallest $t_0$ such that for all $t \in [0, t_{N,\text{end}}]$, the minimum distance between the leading and trailing aircraft is at least $S$, and at $t = t_{N,\text{end}}$, the distance is at least $S_{\text{thr}}$. Accordingly, the minimum $t_0$ is attained either when the leading and trailing aircraft are closest on the common path or when the leading aircraft crosses the runway threshold. The procedure for its computation is illustrated in Fig. 4.

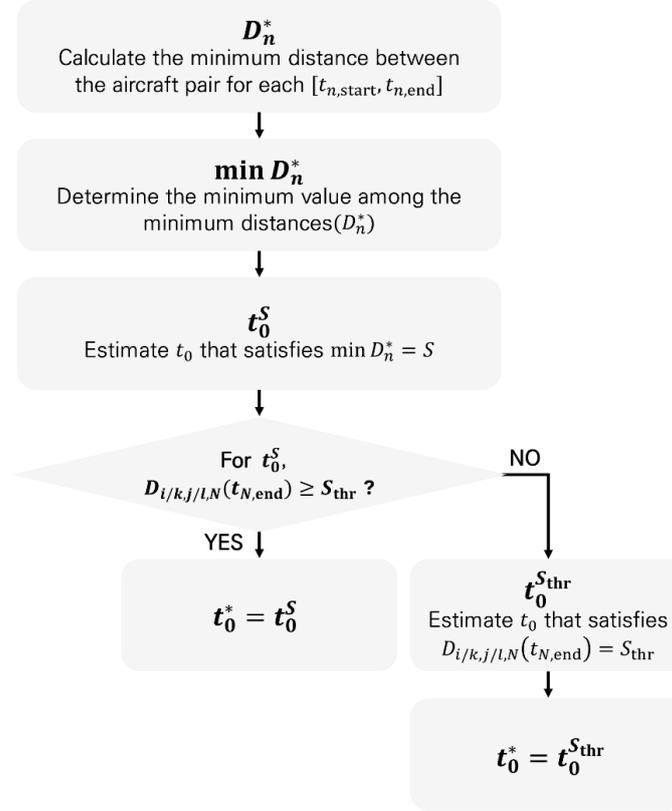

Fig. 4 Flowchart for computing the minimum $t_0$

For each subinterval $[t_{n,\text{start}}, t_{n,\text{end}}]$, the time at which the distance between the leading and trailing aircraft is minimized is determined, and then, the corresponding distance ($D_n^*$) is obtained. Next, the minimum value among all $D_n^*$ (min $D_n^*$) is computed. Subsequently, $t_0$ such that min $D_n^*$ equals $S$ is determined. With this value ($t_0^S$) applied, the minimum $t_0$ ($t_0^*$) is $t_0^S$ when the distance between the leading and trailing aircraft at $t = t_{N,\text{end}}$ is at least $S_{\text{thr}}$. Otherwise, the minimum $t_0$ ($t_0^*$) is $t_0^{S_{\text{thr}}}$, which is the value of $t_0$ that makes the distance at $t = t_{N,\text{end}}$ equal to $S_{\text{thr}}$; in this case, min $D_n^*$ is guaranteed to be no less than $S$.

## 4. Model Application and Analysis

The proposed model was applied to the TMA of Jeju International Airport to estimate its capacity. To this end, the necessary data for the model were collected and preprocessed, and subsequently, they were applied to the model to derive the results. Based on these results, a sensitivity analysis was conducted with respect to the temporal flight distance ($D_{\text{temp}}$), average time separation at the runway threshold ($\bar{T}_{\text{thr}}$), and TMA capacity ($\lambda_{\text{rwy}}$).

4.1 Model Application

Jeju International Airport handles the largest volume of domestic air traffic in Korea. In terms of total air traffic including international operations, it ranks second after Incheon International Airport (Airportal, 2025). The Jeju TMA and Seoul TMA are classified as airports with highly congested TMAs. However, compared with the Seoul TMA, the Jeju TMA is smaller in terms of physical space, and therefore, the structural constraints on handling concentrated traffic demand are correspondingly more pronounced (Ministry of Land, Infrastructure and Transport Republic of Korea, 2023a). In addition, among the two intersecting runways of the airport, RWY 13/31 is used exclusively for departures, whereas landings rely on RWY 07/25, which further impose constraints from a runway-operations perspective (Ministry of Land, Infrastructure and Transport Republic of Korea, 2023b; Lee et al., 2025). Considering these structural and operational constraints, this study selected the RWY 07/25 TMA of Jeju International Airport as the target TMA for model application.

To this end, the arrival flight path was defined based on the most frequently used STARs and IAPs connected to RWY07/25. For RWY07, the Papa STAR and RNP Y IAP were used, whereas Tango STAR and RNP IAP were used for RWY25 (Ministry of Land, Infrastructure and Transport Republic of Korea, 2023b). Accordingly, the arrival flight paths linking the TMA entry points (DOTOL, UPGOS, SOSDO, LIMDI, TAMNA, and TOSAN) to the runway threshold were established, as shown in Figure 5.

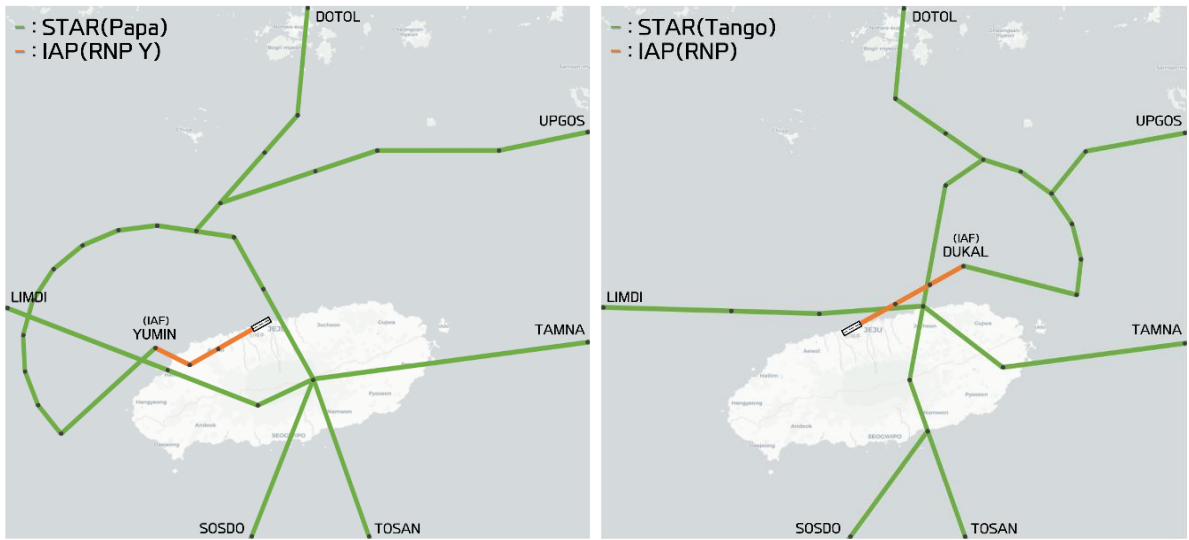
Fig. 5 Arrival flight path of the TMA (left: RWY 07, right: RWY 25)

FlightRadar24 ADS-B trajectory data collected over a one-year period from January to December 2023 were utilized to incorporate the air traffic characteristics of the TMA. The data were preprocessed, and the traffic and aircraft mix proportions for each arrival flight path were extracted (Tables 1 and 2).

Table 1. Proportions for RWY07 TMA

| Entry point | Ratio (%) | | |
| --- | --- | --- | --- |
| | Flight path | A/C type | |
| | | Heavy | Medium |
| DOTOL | 72 | 7 | 93 |
| UPGOS | 22 | - | 100 |
| SOSDO | 4 | 1 | 99 |
| LIMDI | 1 | - | 100 |
| TAMNA | 1 | - | 100 |
| TOSAN | - | - | - |

Table 2. Proportions for RWY25 TMA

| Entry point | Ratio (%) | | |
| --- | --- | --- | --- |
| | Flight path | A/C type | |
| | | Heavy | Medium |
| DOTOL | 72 | 7 | 93 |
| UPGOS | 22 | 2 | 98 |
| SOSDO | 4 | 2 | 98 |
| LIMDI | 1 | - | 100 |
| TAMNA | 1 | - | 100 |
| TOSAN | - | - | - |

The average passing speeds by aircraft type at the entry point, IAP merging point, and runway threshold were extracted for each arrival flight path (Tables 3 and 4). The ground speed recorded in the ADS-B data was used as the aircraft speed to ensure that the model reflects realistic flight times.

Table 3. Average passing speeds for RWY07 TMA (Entry, $MP_{iap}$, and THR)

| Entry point | A/C type | Speed (kt) | | |
|---|---|---|---|---|
| | | Entry | $MP_{iap}$ | THR |
| DOTOL | H | 352 | 196 | 132 |
| | M | 343 | 195 | 143 |
| UPGOS | H | - | - | - |
| | M | 316 | 192 | 140 |
| SOSDO | H | 348 | 206 | 125 |
| | M | 350 | 193 | 135 |
| LIMDI | H | - | - | - |
| | M | 318 | 189 | 137 |
| TAMNA | H | - | - | - |
| | M | 338 | 187 | 143 |
| TOSAN | H | - | - | - |
| | M | - | - | - |

Table 4. Average passing speeds for RWY25 TMA (Entry, $MP_{iap}$, and THR)

| Entry point | A/C type | Speed (kt) | | |
|---|---|---|---|---|
| | | Entry | $MP_{iap}$ | THR |
| DOTOL | H | 311 | 191 | 132 |
| | M | 306 | 189 | 143 |
| UPGOS | H | 295 | 185 | 132 |
| | M | 284 | 187 | 142 |
| SOSDO | H | 396 | 214 | 126 |
| | M | 367 | 191 | 134 |
| LIMDI | H | - | - | - |
| | M | 365 | 192 | 132 |
| TAMNA | H | - | - | - |
| | M | 329 | 181 | 148 |
| TOSAN | H | - | - | - |
| | M | - | - | - |

Table 5 presents the results obtained by applying the extracted data values as input variables to the model. In this process, considering the radar horizontal separation and separation between arrival aircraft applied at Jeju International Airport, the ATC separations for the TMA and runway threshold were set to 5 and 8 NM, respectively (Ministry of Land, Infrastructure and Transport Republic of Korea,

2022; Lee et al., 2025).

Table 5. Results of TMA capacity estimation

| RWY | $D_{\text{temp}}$(min) | $\bar{T}_{\text{thr}}$(min) | $\lambda_{\text{rwy}}$(A/C count) |
|---|---|---|---|
| RWY07 | 28.51 | 3.06 | 9.3 |
| RWY25 | 21.47 | 3.10 | 6.9 |

The TMA capacities of RWY 07 and RWY 25 at Jeju International Airport were estimated to be 9.3 and 6.9 in terms of occupancy count, respectively. The difference of 2.4 arises because the RWY 25 TMA has a larger average time separation at the runway threshold but a shorter temporal flight distance than that of the RWY 07 TMA. The average time separation at the runway threshold differs by 0.04 min, whereas the temporal flight distance differs by 7.04 min, which suggests that the temporal flight distance has a greater effect on capacity difference. These results indicate that, when RWY 07 is used as the arrival runway, the structural space of the Jeju International Airport TMA is formed more extensively than that when RWY 25 is used.

4.2 Sensitivity Analysis

The TMA capacity is determined by variables related to the structure of the airspace and characteristics of air traffic. A sensitivity analysis is essential for evaluating the effect of these variables on TMA capacity. Based on an application case of the Jeju International Airport TMA, a sensitivity analysis was conducted on key variables affecting TMA capacity, temporal flight distance, and average time separation at the runway threshold. In this process, variables difficult to adjust in the short term within the target TMA such as traffic proportions and physical distance by arrival flight path were excluded from the sensitivity analysis. Instead, effects of variables such as the average passing speeds by aircraft type and ATC separations applied within the TMA and runway threshold were analyzed in relation to capacity estimation.

A variation in the range of –10% to +10% was applied based on the extracted average passing speed for the average passing speed by aircraft type. However, such a variation was not applied to the average passing speed at the runway threshold considering flight safety during the landing phase.

The ATC separations applied within the TMA and at the runway threshold were varied considering the radar horizontal separation minima prescribed by ATC procedures (3 and 5 NM) and the separation between arriving aircraft (8 NM) under mixed-use runway operations handling both arrivals and departures at Jeju International Airport (Ministry of Land, Infrastructure and Transport Republic of Korea, 2022; Lee et al., 2025) as follows: (a) $S$ = 5 NM, $S_{thr}$ = 8 NM; (b) $S$ = 5 NM, $S_{thr}$ = 5 NM; (c) $S$ = 3 NM, $S_{thr}$ = 5 NM; and (d) $S$ = 3 NM, $S_{thr}$ = 3 NM.

The results of the sensitivity analysis of temporal flight distance ($D_{temp}$), average time separation at runway threshold ($\bar{T}_{thr}$), and TMA capacity ($\lambda_{rwy}$) with respect to variations in specified variables are presented in Figures 6, 7, and 8.

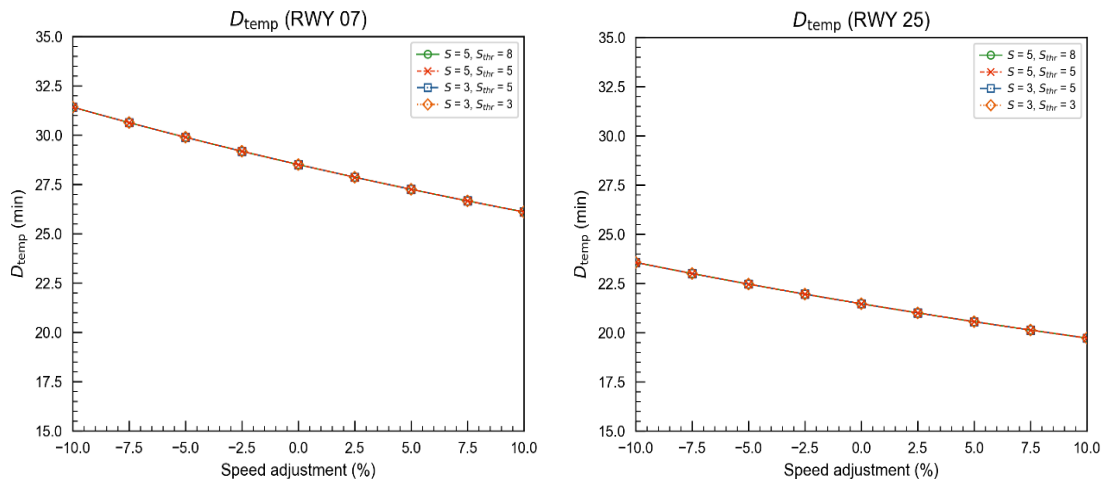

Fig. 6 Sensitivity analysis of the RWY07/25 TMA: $D_{temp}$

Variations in $S$ and $S_{thr}$ in both the RWY07 and RWY25 TMA did not affect $D_{temp}$, whereas it decreased with an increase in the average passing speed by aircraft type. In the RWY07 TMA, the value decreased by ~20.4%, from 31.42 min at a −10% variation to 26.1 min at a +10% variation from the baseline speed. In the RWY25 TMA, it decreased by ~19.5%, from 23.57 min at a −10% variation to 19.72 min at a +10% variation from the baseline speed, indicating a relatively smaller change.

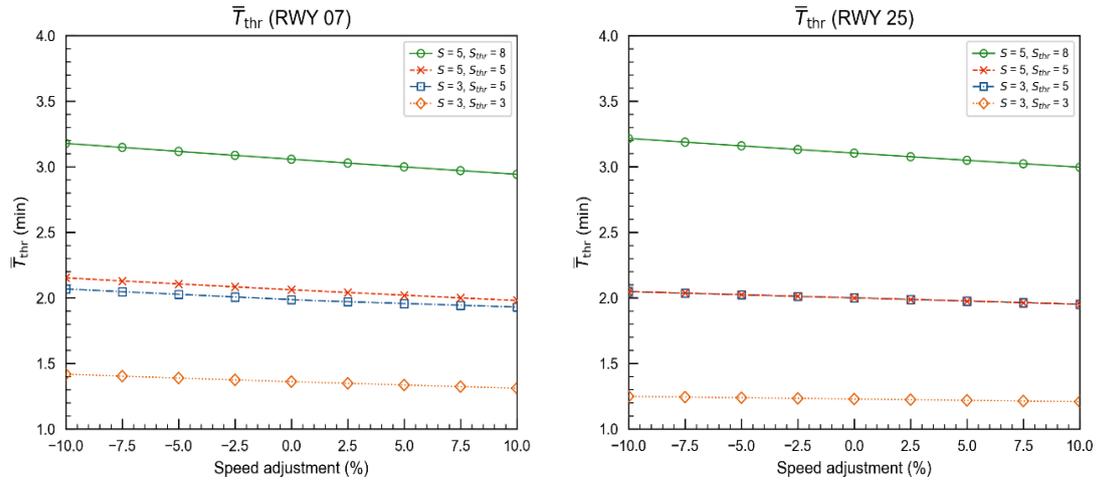

Fig. 7 Sensitivity analysis of the RWY07/25 TMA: $\bar{T}_{thr}$

Further, $\bar{T}_{thr}$ decreased in RWY07 and RWY25 TMA with decreasing $S$ and $S_{thr}$ and increasing average passing speed by aircraft type. Moreover, the slope of $\bar{T}_{thr}$ with respect to variations in speed became steeper with increasing applied ATC separations. $\bar{T}_{thr}$ exhibited the steepest variation at $S = 5$ NM and $S_{thr} = 8$ NM, and the most gradual variation at $S = 3$ NM and $S_{thr} = 3$ NM. For $S = 5$ NM and $S_{thr} = 8$ NM, the RWY 07 TMA decreased by 0.24 min (from 3.18 to 2.94 min), which corresponds to an average slope of 0.012 min per 1% speed change. RWY 25 TMA decreased by 0.22 min (from 3.22 to 3 min) with a slope of 0.011 min per 1% speed change. For $S = 3$ NM and $S_{thr} = 3$ NM, the RWY 07 TMA decreased by 0.11 min (from 1.42 to 1.31 min), which corresponds to an average slope of 0.005 min per 1% speed change, while RWY 25 TMA decreased by 0.04 min (from 1.25 to 1.21 min), which corresponds to an average slope of 0.002 min per 1% speed change.

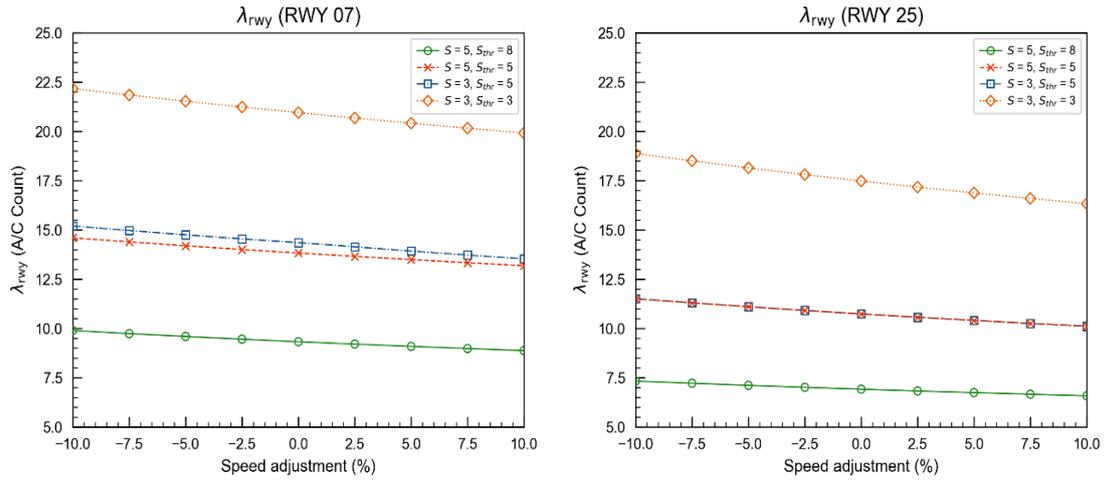

Fig. 8 Sensitivity analysis of the RWY07/25 TMA: $\lambda_{rwy}$

For $\lambda_{rwy}$, both the RWY07 and RWY25 TMA showed an increasing tendency with decreasing $S$ and $S_{thr}$ and decreasing average passing speed by aircraft type. Focusing on the variations in $S$ and $S_{thr}$, the capacity increased because $\bar{T}_{thr}$ decreased with larger applied ATC separations, whereas $D_{temp}$ remained unaffected. In terms of average passing speed, the reduction in $\bar{T}_{thr}$ associated with increasing speed was smaller than the reduction in $D_{temp}$. Therefore, the capacity increased with a decrease in passing speed, which indicates that the effect of reduction in $D_{temp}$ had a relatively greater effect on capacity.

The results of the sensitivity analysis indicated that the TMA capacity is affected by the combined effects of aircraft speed and applied ATC separations. Therefore, in scenarios where it is difficult to make short-term adjustments to traffic characteristics, such as in the traffic proportions for each arrival flight path, or to instrument flight procedures, capacity enhancement may be pursued through operational measures such as speed restrictions or by reducing ATC separations.

## 5. Conclusion

This study analyzed the relationship between the structure of terminal airspace and its capacity and proposed a capacity model reflecting this relationship, which provided a theoretical foundation for the effective utilization of future ATM systems and development of strategies for airspace management. To establish this model, the fundamental concept of TMA capacity was revisited. TMA capacity can be defined either as the maximum number of aircraft passing a specified point in the airspace or as the maximum number of aircraft occupying the airspace. These conceptual definitions differ only in their intended purpose and approach and are best regarded as complementary. This study focused on the structural space of TMA, and under conditions that did not require ATCo intervention (e.g., radar vectoring), defined capacity as the maximum number of aircraft occupying the airspace when aircraft fly the STARs and IAPs. Further, the structural space required for the model was defined conceptually and traffic characteristics were parameterized. Subsequently, they were implemented through mathematical methods to derive the temporal distance and average separation. Finally, a capacity estimation model for TMA was developed by integrating the temporal flight distance and average time separation. The developed model was applied to the TMA of Jeju International Airport, and a sensitivity analysis was conducted to assess the effect of key variables on TMA capacity. The results confirmed that variations in the average passing speeds by aircraft type and applied ATC separations significantly affected capacity. These results indicate that capacity enhancement can be achieved through operational adjustments such as speed restrictions or ATC separation modifications even under conditions where structural constraints within a physically limited TMA cannot be directly alleviated. This contributes to the provision of practical measures for increasing TMA capacity in future operations.

This study recognized TMA as an independent airspace unit rather than as a general airspace sector and defined its capacity based on the concept of maximum occupancy count reflecting both the structure and traffic characteristics of airspace. Furthermore, this study presents a model for the quantitative estimation of this capacity. A limitation of this study is that factors related to ATCo intervention were not considered in the model. Accordingly, future research should aim to improve the model to incorporate ATCo workload within the structural space so that capacity can be estimated considering this factor.